\documentclass{amsart}
\usepackage{amsmath,amssymb,hyperref,mathrsfs,graphicx}
\usepackage[utf8x]{inputenc}
\usepackage{amsmath}
\usepackage{amsfonts}
\usepackage{latexsym}
\usepackage{amsthm}
\usepackage{amssymb,amscd}
\usepackage{xargs}
\usepackage{tikz}
\usepackage{graphicx}
\usepackage{color} 
\usepackage{enumitem}

\newtheorem{thm}{Theorem}[section]

\newtheorem{lem}[thm]{Lemma}
\newtheorem{defn}[thm]{Definition}

\newtheorem{cor}[thm]{Corollary}
\newtheorem{rmk}[thm]{Remark}
\newcommand{\be}{\begin{eqnarray}}
\newcommand{\ee}{\end{eqnarray}}
\newcommand{\bee}{\begin{eqnarray*}}
\newcommand{\eee}{\end{eqnarray*}}
\newcommand{\beal}{\begin{aligned}}
\newcommand{\enal}{\end{aligned}}

\newcommand{\T}{\mathbb{T}}
\newcommand{\R}{\mathbb{R}}
\newcommand{\Q}{\mathbb{Q}}
\newcommand{\cG}{\mathcal{G}}

\newcommand{\N}{\mathbb{N}}
\newcommand{\Z}{\mathbb{Z}}

\newcommand{\cO}{\mathcal{O}}
\newcommand{\cM}{\mathcal{M}}
\newcommand{\cA}{\mathcal{A}}

\newcommand{\cC}{\mathcal{C}}

\newcommand{\cP}{\mathcal{P}}
\newcommand{\cU}{\mathcal{U}}
\newcommand{\cN}{\mathcal{N}}

\newcommand{\cL}{\mathcal{L}}

\title{Generically Ma\~{n}\'e set supports uniquely ergodic measure for residual cohomology class}
\author{Jianlu Zhang}
\address{Department of Mathematics, University of Toronto\\Ontario, Canada}
\email{jianlu.zhang@utoronto.ca}
\thanks{}
\date{}
\begin{document}
\maketitle
\begin{abstract} In this paper, we proved that for generic Tonelli Lagrangian,
there always exists a residual set $\mathcal{G}\subset H^1(M,\mathbb{R})$ such that
\[
\widetilde{\mathcal{M}}(c)=\widetilde{\mathcal{A}}(c)=\widetilde{\mathcal{N}}(c),\quad \forall c\in\mathcal{G}
\]
with $\widetilde{\mathcal{M}}(c)$ supports on a uniquely ergodic measure.
\end{abstract}
\vspace{10pt}

\noindent For positively definite Hamiltonian systems, the Mather theory now proves to be a prominent tool in exploring the dynamic behavior of invariant sets. Its essence is using the variational method to classify a list of invariant sets with different action minimization properties, and these variational properties usually can be transferred into fine topological features. Benefit from this we can construct interesting orbits and reveal the corresponding dynamic phenomena, see \cite{B,Mat,Mat2}. \\

Although this theory was initially discovered by J. Mather in 1980s and used to solve a list of monotone twist map problems, soon more mathematicians applied it to more interesting topics, e.g.   the Arnold Diffusion, the Hamilton-Jacobi equation, the Optimal transportation and etc,  see \cite{CY1,CY2,Fa2,JZ}. Specially mention that R. Ma\~{n}\'e independently developed an action potential approach and pushed this theory forward greatly, see \cite{Mn,Mn2}. He also proposed a list of enlighting conjectures in \cite{Mn2} which inspire this work.\\

This article is organized as follows: In Section \ref{1} we formalized the fundamental constructions of the Mather theory then state our main conclusion. In Section \ref{2} we give the proof and exhibit several applications and corollaries.

\section{\sc Introduction}\label{1}

We consider a compact $n$-dimensional Riemannian manifold $M$ without boundary.  An (autonomous) {\bf Tonelli Lagrangian} on $M$ is a $C^2-$smooth function $L:TM\rightarrow\R$ {with} $(x,v)\in TM$
satisfying these assumptions \cite{Mat}:
\begin{itemize}
\item {\it Convexity:} The Hessian matrix $L_{vv}$ is positively definite for any $(x,v)\in TM$;
\item {\it Superlinearity:} $L(x,v)/\|v\|\rightarrow+\infty$, as $\|v\|\rightarrow+\infty$ for any $x\in M$;
%\item {\bf completeness} the Euler Lagrange equation of $L(x,v,t)$ is well defined for the whole time $t\in\R$;
\end{itemize}
Based on these, the action function
\be
A(x,y,t):=\inf_{\substack{\gamma(0)=x\\
\gamma(t)=y}}\int_0^tL(\gamma(t),\dot{\gamma}(t))dt,\quad t>0
\ee
is well defined for the absolutely continuous curves $\gamma\in\cC^{ac}([0,t],M)$, in other words, the existence of the minimum for a Tonelli Lagrangian can be ensured, see \cite{Mat}. The extremals satisfy the {\bf Euler-Lagrangian equation} which in local coordinates is given by
\be\label{EL}
\frac d {dt}L_v(\gamma(s),\dot{\gamma}(s))=L_x(\gamma(s),\dot{\gamma}(s)), \quad\forall s\in[0,t].
\ee
Recall that we can get the {\bf conjugated Hamiltonian} $H(x,p)$ for $(x,p)\in T^*M$ by the Legendre transformation:
\[
H(x,p)=\max_{v\in T_xM}\{\langle p,v\rangle-L(x,v)\}.
\] 
Since $L(x,v)$ is autonomous, $H(x,p)$ is a natural first integral of the Euler-Lagrangian flow $\phi_L^t$. Moreover, the Euler-Lagrangian flow can be expanded for $t\in\R$ due to the convexity of $H(x,p)$.\\

As now the completeness of the E-L flow holds, we naturally get a set of all the flow-invariant probability measures on $TM$, which can be denoted by $\mathfrak{M}_L$. This is due to the Birkhoff Ergodic Theorem. Then we can classify $\mathfrak{M}_L$ by the cohomology class $c\in H^1(M,\R)$, and define the $\alpha(c): H^1(M,\R)\rightarrow\R$ by
\be\label{alpha}
\alpha(c)=-\inf_{\mu\in\mathfrak{M}_L}\int L-\eta\; d\mu, \quad[\eta]=c.
\ee
Based on the same conjugate principle, we can get $\beta(h): H_1(M,\R)\rightarrow\R$ by
\be
\beta(h)=\inf_{\mu\in\mathfrak{M}_L, \rho(\mu)=h}\int L\;d\mu
\ee
where $\rho(\mu)\in H_1(M,\R)$ is called the {\bf rotational vector} and defined by
\[
\langle [\lambda], \rho(\mu)\rangle=\int\lambda\;d\mu,\quad\forall \text{\;closed 1-form \;}\lambda \text{\;on\;} M.
\]
Due to the positive definiteness and super linearity, both of these two functions are convex and superlinear, and
\[
\langle c,h\rangle\leq\alpha(c)+\beta(h),\quad\forall c\in H^1(M,\R),\; h\in H_1(M,\R),
\]
where the equality holds only for $c\in D^+\beta(h)$ and $h\in D^+\alpha(c)$ (sub-derivative set). We denote by $\mathcal{P}_L(c)\subset\mathfrak{M}_L$ the $c-$minimal measure set and $\widetilde{\mathcal{M}}(c)=$ supp$ \mathcal{P}_L(c)\subset TM$ the closure of the union for all the supports of the minimizng measures of (\ref{alpha}), which is the so called {\bf Mather set}. Its projection to $M$ is the {\bf projected Mather set} $\mathcal{M}(c)$. 
From \cite{Mat} we know that $\pi^{-1}\big{|}_{\mathcal{M}(c)}:M\rightarrow TM$ is a Lipschitz graph, where $\pi$ is the standard projection. 
\begin{rmk}
Based on R. Ma\~{n}\'e's setting in \cite{Mn}, we can enlarge the variational space from $\mathfrak{M}_L$ to  the set of closed probability measures $\mathfrak{M}_{cl}$. Each $ \mu_{cl}\in\mathfrak{M}_{cl}$ can be uniquely decided due to the Birkhoff ergodic theorem:
\[
\int f d\mu_{cl}:=\frac{1}{T_{cl}}\int_0^{T_{cl}}f(\gamma_{cl},\dot{\gamma}_{cl})\;dt,\quad\forall f\in C^{ac}(TM,\R)
\]
where $T_{cl}$ is the period of the loop $\gamma_{cl}$, but $\dot\gamma^-(T_c)\neq\dot\gamma^+(0)$ may happen. We can still get the same $\alpha(c)$ and $\beta(h)$ under this new setting.
\end{rmk}

To make a preciser portrait of the phase space, R. Ma\~{n}\'e defined the so-called {\bf Action Potential function}, which is shown as
\begin{equation}
\Phi_c(x,y)=\inf_{t\geq0}h_c(x,y,t)
\end{equation}
with
\begin{equation}
h_c(x,y,t)=\inf_{\substack{\xi\in C^{ac}([0,t],M)\\
\xi(0)=x\\
\xi(t)=y}}A_c(\xi)\big{|}_{[0,t]},
\end{equation}
and
\begin{equation}
A_c(\gamma)\big{|}_{[0,t]}=\int_0^{t}L(\gamma(s),\dot{\gamma}(s))-\langle\eta_c(\gamma(s)),\dot{\gamma}(s)\rangle ds+\alpha(c)t.
\end{equation}
 Then a curve $\gamma:\mathbb{R}\rightarrow M$ is called {\bf c-semi static} if 
\[
\Phi_c(\gamma(a),\gamma(b))=A_c(\gamma)\big{|}_{[a,b]},
\]
for all $a,b\in\mathbb{R}$. A semi static curve $\gamma$ is called {\bf c-static} if
\[
A_c(\gamma)\big{|}_{[a,b]}+\Phi_c(\gamma(b),\gamma(a))=0,\quad\forall a,b\in\mathbb{R}.
\]
The {\bf Ma\~{n}\'e set} $\widetilde{\mathcal{N}}(c)\subset TM$ is denoted by the set of all the c-semi static orbits, and the {\bf Aubry set} $\tilde{\mathcal{A}}(c)$ is the set of all the c-static orbits. From \cite{B} we can see that $\forall x,y\in M$, 
\be\label{class}
d_c(x,y):=\Phi_c(x,y)+\Phi_c(y,x)\geq0
\ee
always holds, which implies that every static curve should be a semi-static curve first. Moreover, $\pi^{-1}:\mathcal{A}(c)\rightarrow \tilde{\mathcal{A}}(c)$ is also a Lipschitz graph, but usually $\widetilde{\mathcal{N}}(c)$ is not a graph over ${\mathcal{N}}(c)$. \\

Before we explore the further features for the aforementioned variational minimal sets, we introduce the following definition first:
\begin{defn}
We say a property is {\bf generic} if for any fixed Tonelli Lagrangian $L$, there exists a residual (countable intersection of open and dense subsets) set $\cO\subset C^\infty(M,\R)$, such that the property holds for $L + f$, $\forall f\in \cO$.
\end{defn}
\begin{lem}[R. Ma\~{n}\'e, \cite{Mn}]\label{lem1}
For a fixed $c\in H^{1}(M,\R)$, there exist residual subsets $\cO_{c}\subset C^\infty(M, \R)$ such that $\#\mathfrak{M}_{L+\psi}(c) = 1$, i.e. there exists a unique ergodic $c-$minimal measure.

Moreover, for  a fixed $h\in H_{1}(M,\R)$, there exist residual subsets $\cO_{h}\subset C^\infty(M, \R)$ such that $\#\mathfrak{M}_{L+\psi}(c_{h}) = 1$, i.e. there exists a unique ergodic $c_{h}-$minimal measure, where $c_{h}\in D^{+}\beta(h)$.
\end{lem}
\begin{lem}[P. Bernard, \cite{B}]\label{lem2}
\begin{itemize}
\item $\widetilde{\mathcal{M}}(c)\subset\tilde{\mathcal{A}}(c)\subset\widetilde{\mathcal{N}}(c)$.
\item If $\#\mathfrak{M}_{L}(c) = 1$, then $\tilde{\mathcal{A}}(c)=\widetilde{\mathcal{N}}(c)$.
\end{itemize}
\end{lem}
\begin{figure}
\begin{center}
\includegraphics[width=8cm]{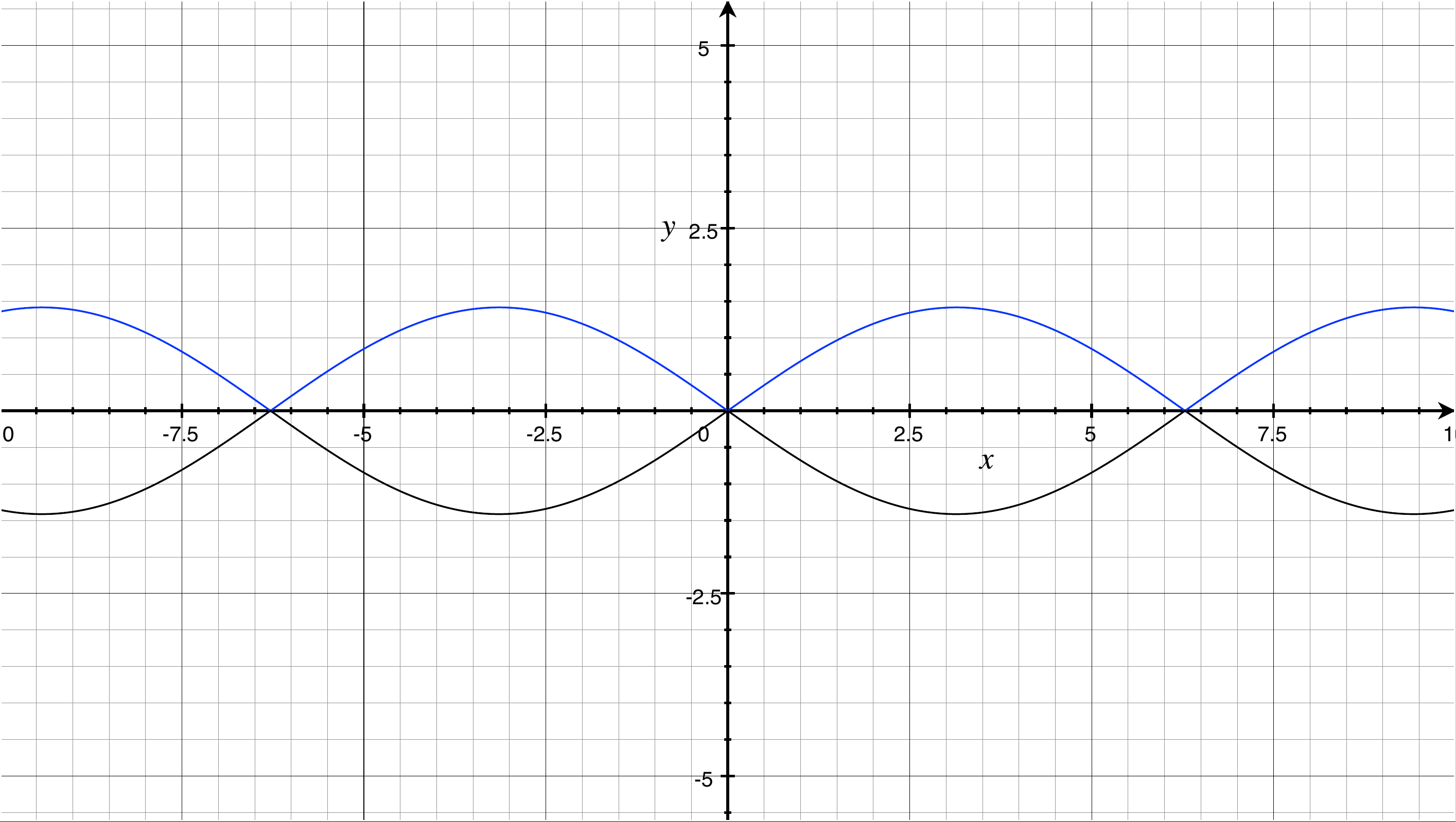}
\end{center}
\begin{center}
\includegraphics[width=8cm]{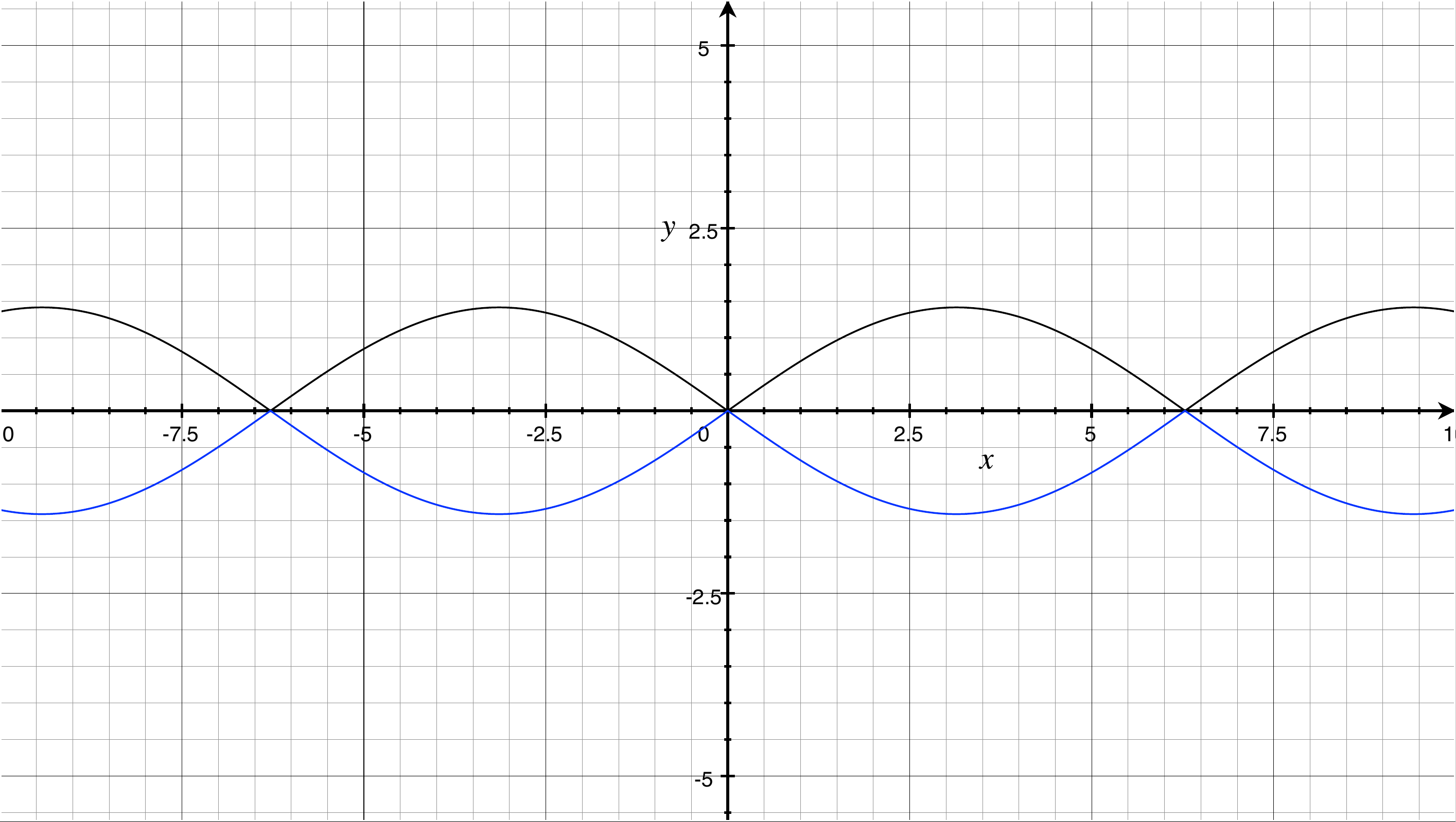}
\caption{For $c=4/\pi$, the Aubry set consists of the hyperbolic fixed point and the upper separatrix. For $c=-4/\pi$, the Aubry set consists of the hyperbolic fixed point and the lower separatrix.}

\label{fig}
\end{center}
\end{figure}

\begin{rmk}
Usually $\widetilde{\mathcal{M}}(c)$ may not equal $\widetilde{\mathcal{A}}(c)$. For the pendulum $H(x,y)={1}/{y^2}+(\cos x-1)$, $(x,y)\in T^*\T$, we can see that $\alpha(c)$ is $C^1$ smooth and $\alpha(c)\geq0$ with `$=$' holds for $c\in[-4/\pi,4\pi]\subset H^1(\T,\R)$. But $\widetilde{\mathcal{M}}(4/\pi)=\{(0,0)\}$ is strictly contained in $\tilde{\mathcal{A}}(4/\pi)=\{(x,\sqrt{2(1-\cos x)})|x\in\T\}$ and $\widetilde{\mathcal{M}}(-4/\pi)=\{(0,0)\}$ is strictly contained in $\tilde{\mathcal{A}}(-4/\pi)=\{(x,-\sqrt{2(1-\cos x)})|x\in\T\}$ (see Figure \ref{fig}).

Another counter example is that for the irrational $h\in H_1(M,\R)$, $\widetilde{\mathcal{M}}(c_h)\subsetneq\tilde{\mathcal{A}}(c_h)$ could still happen for $c_h\in D^+\beta(h)$ (Appartently there are uncountablely many $h$ could be chosen). This point is exposed by M. Arnaud in \cite{Ar}, which showed that for a fixed $w$ irrational number, there exists an exact monotone twist map such that $\widetilde{\cM}(c_w)$ is a Denjoy set but $\widetilde{\cA}(c_w)$ is an invariant circle. 
\end{rmk}
\vspace{10pt}

With all these evidences above, a natural question arises: To what degree can we read out the dynamic informations of $\widetilde{\cN}(c)$ only from $\widetilde{\cM}(c)$ ? To answer this, we state our main result as following:
% \bthm\label{MR}
%\ethm
% \includegraphics[height=10cm]{gsf.png}
\begin{thm}[Main Conclusion]
For generic Tonelli Lagrangian $L(x,v)$ defined on $TM$, there always exists a residue set $\mathcal{G}^*\subset H^1(M,\mathbb{R})$ such that
\[
\widetilde{\mathcal{M}}(c)=\widetilde{\mathcal{A}}(c)=\widetilde{\mathcal{N}}(c),\quad \forall c\in\mathcal{G}^*
\]
with $\widetilde{\mathcal{M}}(c)$ supports a uniquely ergodic measure.
\end{thm}
From the viewpoint of topology dynamics, this result greatly reduces the complexity of the Ma\~{n}\'e set, at least for a `big' type of cohomology classes.
\begin{defn}
We call an invariant set $\Omega$ be {\bf topological minimal} if $\Omega$ could not be further decomposed into a union of smaller invariant sets. 

%An invariant set $\Omega$ is called {\bf chain transitive} if for any two $x,y\in\Omega$ and any  $\epsilon > 0$, $T > 0$ there exists a {\bf $(\epsilon, T)$-chain} in $\Omega$ joining $x$ and $y$, in the sense that we can find a finite sequence $\{(z_i,t_i)\}_{i=1}^N$ such that $x=z_1$, $y=z_N$, $t_i>T$ and 
%$d(\phi_H^{t_i}(z_i),z_{i+1})<\epsilon$ for $i=1,\dots, N-1$. When this condition holds only for $x=y\in\Omega$ we say that $\Omega$ is {\bf chain recurrent}.
\end{defn}
\begin{cor}
For generic Tonelli Lagrangian, there always exists a residual set $\mathcal{G}^*\subset H^1(M,\mathbb{R})$ such that $\widetilde{\mathcal{N}}(c)$ is a topological minimal set which is a Lipschitz graph as well.
\end{cor}
However, we have to confess that although $\mathcal{G}$ is a topologically `big' set of which the Ma\~{n}\'e set becomes rather `regular', its geometric structure could be very complicated. As an enlightening supplement, we will talk about this point later on. 

%From current situation \cite{Mn,Mn2}, hyperbolicity seems to be useful to reduce the complexity of the topological structure of $\mathcal{G}$. We will talk about this point later on.
\section{\sc Proof of the theorem}\label{2}

Let's fix the Tonelli Lagrangian by $L_{0}$ in this part. Due to Lemma \ref{lem1}, we can take a sequence $\{c_n\}_{n=1}^\infty\subseteq H^1(M,\mathbb{R})$ which is dense in $H^1(M,\mathbb{R})$, such that there exists a residual set $\mathcal{O}^{'}\subset C^{\infty}(M,\R)$ being the perturbation of $L_{0}$, such that $\mathcal{P}(L_{0},\psi,c_n)$ is uniquely ergodic for all $\{c_{n}\}_{n=1}^{\infty}$ and $\psi\in\mathcal{O}'$.
\vspace{15pt}

{\bf Step 1:} For a fixed $c_n$ and $ \psi\in \mathcal{O}^{'}$, from Lemma \ref{lem2} we get $\widetilde{\mathcal{A}}(L_{0},\psi,c_n)=\widetilde{\mathcal{N}}(L_{0},\psi,c_n)$. We claim that for a sufficiently small $\varphi\in C^{\infty}(M,\mathbb{R})$,
\[
\widetilde{\mathcal{M}}(L_{0},\psi+\varphi,c_n)=\widetilde{\mathcal{A}}(L_{0},\psi+\varphi,c_n)=\widetilde{\mathcal{N}}(L_{0},\psi+\varphi,c_n).
\]
This is easily achievable because we just need to take $\varphi:M\rightarrow\R$ by
\begin{equation*}
\varphi(x)=\begin{cases}
0,\quad \text{on }\ {\mathcal{M}}(L_{0},\psi,c_n),\\
>0\text{\;and\;}\ll 1,\    \text{on }\ {\mathcal{M}}(L_{0},\psi,c_n)^c.
%\ ,\text{where}\ \mathcal{U}\ \text{is a sufficient small neighborhood of}\ {\mathcal{M}}_L(c_n).\\
\end{cases}
\end{equation*}
Based on these, $\widetilde{\mathcal{M}}(L_{0},\psi+\phi,c_n)=\widetilde{\mathcal{M}}(L_{0},\psi,c_n)$. Recall that the Graph Property holds for $\widetilde{\mathcal{A}}(L_{0},\psi+\phi,c_n)$, if ${\mathcal{A}}(L_{0},\psi+\phi,c_n)\backslash{\mathcal{M}}(L_{0},\psi+\phi,c_n)\neq\emptyset$,
at least one point $x'$ exists such that 
\be
\Phi_{c_n}^{L_{0}+\psi+\varphi}(y',x')+A_{c_n}^{L_{0}+\psi+\varphi}(\gamma')|_{[0,1]}&\geq&\Phi_{c_n}^{L_{0}+\psi}(y',x')+\int_0^1\varphi(\gamma'(s))ds+\Phi_{c_n}^{L_{0}+\psi}(x',y')\nonumber\\
&>&\Phi_{c_n}^{L_{0}+\psi}(y',x')+\Phi_{c_n}^{L_{0}+\psi}(x',y')\nonumber\\
&\geq&0
\ee
where $\gamma':\R\rightarrow M$ is the Lagrangian flow with $\gamma'(0)=x'$ and $\gamma'(1)=y'$. The first inequality holds as $\varphi(x)$ is non-negative. The second inequality is due to the strict positiveness of $\varphi$ on ${\mathcal{M}}_L(c_n)^c$, as $x'\in{\mathcal{A}}(L_{0},\psi+\phi,c_n)\backslash{\mathcal{M}}(L_{0},\psi+\phi,c_n)$ and $\int_0^1\varphi(\gamma'(s))ds>0$. The last inequality can be derived from (\ref{class}).\\

Recall that $\varphi$ can be made sufficiently small, so we can get a conclusion that:
\vspace{5pt}

{\it for every $c_n$, there exists a dense set $\mathcal{O}_n\subset C^{\infty}(M,\R)$ as the perturbation of $L_{0}$, such that
\be\label{target}
\widetilde{\mathcal{M}}(L_{0},\phi_{n},c_n)=\widetilde{\mathcal{A}}(L_{0},\phi_{n},c_n)=\widetilde{\mathcal{N}}(L_{0},\phi_{n},c_n)
\ee
and $\mathcal{P}(L_{0},\phi_{n},c_{n})$ is uniquely ergodic for all $ \phi_{n}\in\mathcal{O}_{n}$.
}
\vspace{10pt}

{\bf Step 2:} In this step, we get the convergence of these sets under the Hausdorff distance. Before we doing that, the following useful Lemmas should be involved:
\begin{lem}\cite{CY1,CY2}
As a set-valued function, $(\xi,c)\rightarrow\widetilde{\mathcal{N}}(L,\xi,c)$ is upper-semicontinuous w.r.t $\xi\in C^{\infty}(M,\R)$ and $c\in H^{1}(M,\R)$, where we adopt the $C^\infty$ topology and Euclid norm $|\cdot|_e$ each, and the Hausdorff distance on $TM$ (see (\ref{d-h}) for the definition).
\end{lem}
\begin{rmk}
Notice that usually $\widetilde{\mathcal{N}}(L,\xi,c)$ is not upper-semicontinuous! This point is crucial in construct local connecting orbits in nearly integrable systems, see \cite{CY1,CY2} for more details.
\end{rmk}
\begin{lem}\cite{Mn}
In the same setting as above, $\mathcal{P}(L,\xi,c)$ is upper-semicontinuous of $(\xi,c)$ as well under the weak* topology of invariant measure space.
\end{lem}
As $L_{0}$ is fixed once for all, so we can remove it in the notations for short. Then for any two sequences $c_k^n\longrightarrow c_n$ and $\phi_k^n\longrightarrow \phi^{n}$ as $k\rightarrow\infty$ with $\phi^{n}\in\mathcal{O}_n$, we have:
\be\label{dist}
\overline{\lim}\widetilde{\mathcal{A}}(\phi_k^n, c_k^n)\subseteq \overline{\lim}\widetilde{\mathcal{N}}(\phi_k^n, c_k^n)\subseteq \widetilde{\mathcal{N}}(\phi_n, c_n)=\widetilde{\mathcal{M}}(\phi_n, c_n).
\ee
On the other hand, the weak limit of $\mathcal{P}(\phi_k^n, c_k^n)$ must be in $\mathcal{P}(\phi_n, c_n)$. Recall that there is only one unique ergodic measure $\mu(\phi_{n},c_{n})\in{\mathcal{P}}(\phi_n, c_n)$, so
\[
\mu(\phi_k^n,c_k^n)\rightharpoonup \mu(\phi_n,c_n),\quad\text{as\;}k\rightarrow\infty
\]
for any $\mu(\phi_k^n,c_k^n)\in\mathcal{P}(\phi_k^n, c_k^n)$. Here the weak* convergence implies
\be\label{m-dist}
\lim_{k\longrightarrow\infty}\sup_{z\in\widetilde{\mathcal{M}}(\phi_n, c_n)}d(z,\widetilde{\mathcal{M}}(\phi_k^n, c_k^n))=0.
\ee
We define the {\bf Hausdorff distanse} of two sets by:
\be\label{d-h}
d_H(A,B)=max\Big{\{}\sup_{x\in A}d(x,B),\ \sup_{x\in B}d(x,A)\Big{\}}.
\ee
%From this definition we can see that 
%\[
%d_H(A',B)\leq d_H(A,B),\quad\text{ if\;} A'\subset A.
%\]
Then
\be
\lim_{k\longrightarrow\infty}\sup_{z\in\widetilde{\mathcal{M}}(\phi_k^n, c_k^n)}d(z,\widetilde{\mathcal{M}}(\phi_n, c_n))&\leq&\lim_{k\longrightarrow\infty}\sup_{z\in\widetilde{\mathcal{A}}(\phi_k^n, c_k^n)}d(z,\widetilde{\mathcal{M}}(\phi_n, c_n))\nonumber\\
&\leq&\lim_{k\longrightarrow\infty}\sup_{z\in\widetilde{\mathcal{N}}(\phi_k^n, c_k^n)}d(z,\widetilde{\mathcal{M}}(\phi_n, c_n))\nonumber\\
&=&0,
\ee
which is because the upper semi-continuity of the Ma\~{n}\'{e} set, see (\ref{dist}).
On the other side,
\be
\lim_{k\longrightarrow\infty}\sup_{z\in\widetilde{\mathcal{A}}(\phi_n, c_n)}d(z,\widetilde{\mathcal{A}}(\phi_k^n, c_k^n))\leq\lim_{k\longrightarrow\infty}\sup_{z\in\widetilde{\mathcal{M}}(\phi_n, c_n)}d(z,\widetilde{\mathcal{M}}(\phi_k^n, c_k^n))=0
\ee
due to (\ref{m-dist}).
%Recall that 
%\[
%\widetilde{\mathcal{M}}(L_n,c_n)=\widetilde{\mathcal{A}}(L_n,c_n)=\widetilde{\mathcal{N}}(L_n,c_n)
%\]
%holds and ${\mathcal{P}}(L_n,c_n)$ is unique ergodic, 
So we make the second claim by:
\be\label{conv}
\lim_{k\rightarrow\infty}d_H(\widetilde{\mathcal{A}}(\phi_k^n,c_k^n),\widetilde{\mathcal{A}}(\phi_n,c_n))=\lim_{k\rightarrow\infty}d_H(\widetilde{\mathcal{M}}(\phi_k^n,c_k^n),\widetilde{\mathcal{M}}(\phi_n,c_n))=0.
\ee
In other words, $\widetilde{\mathcal{A}}(\phi,c)$ and $\widetilde{\mathcal{M}}(\phi,c)$ are both continuous in $\phi\in C^{\infty}(M,\R)$ and $c\in H^{1}(M,\R)$ of the Hausdorff distance, as set-valued functions.
\vspace{10pt}

{\bf Step 3: }Now we try to make (\ref{target}) be true for generic potential functions in $C^{\infty}(M,\R)$ and cohomology classes in $H^{1}(M,\R)$. We claim that there exist $\mathcal{O}_{n,r}$ being an open neighborhood of $\mathcal{O}_n$ and $\mathcal{G}_{n,r}$ being an open neighborhood of $c_n$, such that
\[
d_H(\widetilde{\mathcal{A}}(\phi_{n}^{r},c),\widetilde{\mathcal{M}}(\phi_{n}^{r},c))<\frac{1}{r},\quad r\in\mathbb{Z}_{+}
\]
for all $\phi_{n}^{r}\in \mathcal{O}_{n,r}$ and $c\in \mathcal{G}_{{n,r}}$. This is because the convergence of (\ref{conv}) and the triangle inequality of the Hausdorff distance:
\[
d_H(\widetilde{\mathcal{A}}(\phi,c),\widetilde{\mathcal{M}}(\phi,c))\leq d_H(\widetilde{\mathcal{A}}(\phi,c),\widetilde{\mathcal{A}}(\phi_n,c_n))+d_H(\widetilde{\mathcal{M}}(\phi,c),\widetilde{\mathcal{M}}(\phi_n,c_n)).
\]
Then $\mathcal{O}_{n,r}$ is an open-dense set and $\bigcap_{n}{\mathcal{O}}_{n,r}$ is a residue set. Now $\forall c\in\bigcup_n\mathcal{G}_{n,r},\phi\in\bigcap_{n}{\mathcal{O}}_{n,r}$,
\[
d_H(\widetilde{\mathcal{A}}(\phi,c),\widetilde{\mathcal{M}}(\phi,c))<\frac{1}{r}.
\]
Then $\bigcap_{n,r}{\mathcal{O}}_{n,r}$ and $\bigcap_{r}\bigcup_n\mathcal{G}_{n,r}$ both become residue sets and satisfy
\be
\widetilde{\mathcal{A}}(\phi,c)=\widetilde{\mathcal{M}}(\phi,c),\quad \forall \phi\in\bigcap_{\substack{n\in\N\\ r\in\Z_+}}{\mathcal{G}}_{n,r},\ \ c\in\bigcap_{r\in\N}\bigcup_{n\in\N}\mathcal{G}_{n,r}.
\ee
Let's denote $\mathcal{O}''=\bigcap_{n,r}{\mathcal{O}}_{n,r}$ and $\mathcal{G}^*=\bigcap_{r}\bigcup_{n}\mathcal{G}_{n,r}$ for short,  then for $\mathcal{O}^{*}:=\mathcal{O}^{'}\cap\mathcal{O}^{''}$ still being a residue set, such that 
\[
\widetilde{\mathcal{M}}(\phi,c)=\widetilde{\mathcal{A}}(\phi,c)=\widetilde{\mathcal{N}}(\phi,c),\quad \forall \phi\in\mathcal{O}^{*},\;c\in\mathcal{G}^*
\]
holds. Then we finally finish the proof of the main Theorem.
%\[
%\widetilde{\mathcal{M}}(c)=\widetilde{\mathcal{A}}(c),\quad \forall L\in\mathcal{G}^{*}\ and\ c\in \bigcap_{r}\bigcup_n{G}_{n,r}\doteq G^{*}.
%\]
%In [Mn], Ma\~{n}\'{e} constructed a distance of invariant measure space and $c\rightarrow diam\ \mathcal{P}_c$ is a
%upper-semicontinuous map. So the set of continuous points is a residue set. And we recall that
%\[
% diam\mathcal{P}_{c_n}=0,\quad\forall L\in\mathcal{G}^{*},\ and \ c_n\in\{c_n\},
%\]
%so there is a residue set $G_{L}$ of cohomology space such that
%\[
%diam\mathcal{P}_{c}=0,\quad \forall L\in\mathcal{G}^{*}\ and\ c\in G_{L}.
%\]
%So we take $G_{*}\cap G_L$ as a residue set and satisfies our theorem.
\begin{rmk}
For $n=1$, actually we can take $\mathcal{G}^*$ by
\[
\mathcal{G}^*=int \Big{\{}c\in H^1(M,\R)\Big{|}D^+\alpha(c)\in\Q\Big{\}}.
\]
Recall that alpha function is $C^1$ smooth in this case, so $\mathcal{G}^*$ is actually open dense, see \cite{Mat2}.

For $n=2$, the minimal measure can only support on fixed point, periodic orbit, Denjoy set or invariant torus. According these we can classify $\cG^*$ separately:
\begin{itemize}
\item If $c\in\cG^*$ with $\cP(c)$ supporting on a unique periodic orbit with the homology class $h\in H_1(M,\R)$, there must exist a 1-dimensional flat $\mathcal{L}_h\subset \cG^*$, such that 
\[
\langle h,c-c'\rangle=0,\quad\forall c,c'\in\cL_h.
\]
This is because the upper semi-continuity of Ma\~{n}\'e set and 
\bee
\alpha(c')&=&-\int L-c'd\mu_{c'}= \int c'-Ld\mu_{c'}\\
&\geq&\int c'-L d\mu_c\\
&=&-\int L-c d\mu_c+\langle c'-c,h\rangle\\
&=&\alpha(c),
\eee
and we can switch the position of $c$ and $c'$ and finally $\alpha(c')=\alpha(c)$. Recall that the interior of a  flat shares the same Aubry set, then $\mathcal{L}_h\subset \cG^*$ holds, see \cite{Ber}.

We need to specially remark that when the periodic orbit collapses to be a unique fixed point, then the homology class $h=0$ and $\cL_h\subset \cG^*$ becomes 2-dimensional.
\item If $\cP(c)$ supports on an invariant torus, $\alpha'(c)$ is unique and $D^+\beta(\alpha'(c))$ contains only $\{c\}$.
\item If $\cP(c)$ supports on a Denjoy set, $D^+\beta(\alpha'(c))$ may contain a flat. We couldn't exclude this case, which also cause the same difficulty in the so called {\bf Ma\~{n}\'e Conjecture}:\\

{\tt For generic Lagrangian on a closed manifold $M$, there exists an open dense set 
$\cU\in H^{1}(M, \R)$ such that $\forall c\in\cU$, $\widetilde{\cM}(c)$ consists of a single 
periodic orbit, or fixed point.\hspace{300pt}(*)}

Aforementioned evidence shows that, a clearer portrait of the dynamic mechanism of Denjoy sets is a necessary step towards this Conjecture.
\end{itemize}
%In \cite{Mn}, R. Ma\~{n}\'e proposed the following conjecture:\\
%
%{\tt For generic Lagrangian on a closed manifold $M$, there exists a dense open set $\cU\in H^{1}(M, \R)$ such that $\forall c\in\cU$, $\widetilde{\cM}(c)$ consists of a single periodic orbit, or fixed point.\hspace{330pt}(*)}
%
%Till now, this conjecture is still open. A heuristic idea is that for all $(n-1)$-resonant homology classes in $H_{1}(M,\R)$, $\widetilde{\cM}(c)$ consists of a single periodic orbit or fixed point for generic Lagrangian, see Lemma \ref{lem1}. Recall that a vector $v\in\R^{n}$ is called {\bf $d-$resonant} if we can find $d$ linear independent rational vectors $k_{i}\in\Q^{n}$, $i=1,\cdots,d$, such that $\langle k_{i},v\rangle=0$. Then the set of all $(n-1)-$resonant classes, which can be denoted by $\mathcal{H}_{n-1}$ is dense in $H_1(M,\R)$. The conjugated set $\cC_{n-1}$, which is defined by
%\[
%\cC_{n-1}:=\Big{\{}c\in H^1(M,\R)\Big{|}\exists h\in\mathcal{H}_{n-1},\;c\in D^+\beta(h)\Big{\}}.
%\]
%Now we reach the case that {\tt (*)} is true for $\cC_{n-1}$. So we just need to prove that $\cC_{n-1}$ is generic. This point is true for $n=1$, see \cite{Mat2}. For $n\geq2$, 
\end{rmk}

%\section*{Appendix}

\noindent\textbf{Acknowledgement} This article is finally revised in the author's postdoc session and here I show my thanks to the hospitality of the Math. Department of the University of Toronto. I also thank Prof. C-Q. Cheng for several conversations and useful suggestions in this topic.


\begin{thebibliography}{9999}
\bibitem{Ar} Arnaud MC., {\it A non-differentiable essential irrational invariant curve for a $C^1$ symplectic twist map}, Journal of Modern Dynamics, 2011, 5(3): 583-591
\bibitem{B} Bernard P., {\it Connecting orbits of time dependent Lagrangian systems}, Ann. Inst. Fourier, Grenoble 52, 5 (2002), 1533-1568
\bibitem{CY1} Cheng C-Q.\& Yan J., {\it Existence of diffusion orbits in a priori unstable Hamiltonian systems}, J. Differential Geometry , 67 (2004) 457-517.
\bibitem{CY2} Cheng C-Q.\& Yan J., {\it Arnold diffusion in Hamiltonian Systems: a priori unstable case}, J. Differential Geometry, 82 (2009) 229-277.
\bibitem{C} G. Contreras, J. Delgado\& R. Iturriaga, {\it Lagrangian flows II: The dynamics of globally minimizing orbits}. Boletim da Sociedade Brasileira de Matem\'atica 1997, 28(2):155-196
\bibitem{Fa} A. Fathi, {\it Weak KAM Theorem in Lagrangian Dynamics}, Book to appear, the Cambridge Press.
\bibitem{Fa2} A. Fathi\& A. Siconolfi, {\it Existence of $C^1$ critical subsolutions of the Hamilton-Jacobi equation}, Inventiones Mathematicae, 2004, 155(2): 363-388
 \bibitem{Mat} J.N.Mather. {\it Action minimizing invariant measures for postive definite Lagrangian
 systems}. Math. Z. \textbf{207} (1991), 169-207.
 \bibitem{Mat2} J.N. Mather., {\it Variational construction of connecting orbits}, Ann. Inst. Fourier, 43 (1993), 1349- 1368.
\bibitem{Mn} R.Ma\~{n}\'{e}, {\it Generic properties and problems of minimizing measures of Lagrangian systems}. Nonlinearity. \textbf{9} (1996), 273-310.
\bibitem{Mn2} R.Ma\~{n}\'{e}, {\it Lagrangian flows: The dynamics of globally minimizing orbits}, Bol. Soc. Bras. Mat, {\bf28} (1997) 141-153
\bibitem{Ro} R.T.Rockafellar. {\it Convex Analysis}. Princeton University Press, 1970.
\bibitem{JZ} J. Li\& J. Zhang {\it Optimal Transportation for Generalized Lagrangian}, International Journal of Control \& Automation, 2013, 6(1):32-40
\end{thebibliography}
\end{document}